\documentclass[11pt,reqno]{amsart}
\usepackage{amsmath}
\usepackage{amsxtra}
\usepackage{amscd}
\usepackage{amsthm}
\usepackage{amsfonts}
\usepackage{amssymb}
\usepackage{eucal}
\numberwithin{equation}{section}
\newtheorem{thm}{Theorem}[section]

\newtheorem{defn}[thm]{Definition}
\newtheorem{lemma}[thm]{Lemma}
\newtheorem{conj}[thm]{Conjecture}

\newtheorem{remark}[thm]{Remark}
\usepackage{pstricks}
\usepackage{amsmath}
\usepackage{amsxtra}

\usepackage{amscd}
\usepackage{amsthm}
\usepackage{amsfonts}
\usepackage{amssymb}
\usepackage{eucal}

\newcommand{\g}{{\mathfrak{g}}}

\newcommand{\h}{{\mathfrak{h}}}

\newcommand{\n}{{\mathfrak{n}}}
\renewcommand{\sl}{{\mathfrak{sl}}}

\newcommand{\ch}{{\rm ch}}

\newcommand{\ad}{{\rm ad} }
\newcommand{\R}{{\mathbb R}}
\newcommand{\C}{{\mathbb C}}

\newcommand{\Z}{{\mathbb Z}}
\newcommand{\N}{{\mathbb N}}

\newcommand{\qbin}[2]
{{
\left[
\begin{matrix}{\displaystyle #1}\\
{\displaystyle #2}\end{matrix}
\right]
}}

\newcommand{\cA}{{\mathcal A}}
\newcommand{\cC}{{\mathcal C}}

\newcommand{\cF}{{\mathcal F}}
\newcommand{\cH}{{\mathcal H}}

\newcommand{\cU}{{\mathcal U}}
\newcommand{\mult}{{\rm mult}}

\newcommand{\bR}{{\mathbf R}}

\newcommand{\bl}{{\mathbf l}}
\newcommand{\bm}{{\mathbf m}}
\newcommand{\bn}{{\mathbf n}}

\newcommand{\bx}{{\mathbf x}}
\newcommand{\by}{{\mathbf y}}
\newcommand{\bz}{{\mathbf z}}

\newcommand{\hV}{\widehat{V}}
\newcommand{\bzeta}{{\boldsymbol{\zeta}}}

\newcommand{\coinv}{\cC_{\lambda,\bR}[\bm]}

\newcommand{\bmu}{{\boldsymbol{\mu}}}
\newcommand{\bnu}{{\boldsymbol{\nu}}}

\newcommand{\chq}{{\rm ch}_q}

\DeclareMathOperator{\Gr}{Gr}

\usepackage{graphicx}
\newcommand{\KR}{{\rm KR}}
\newcommand{\beq}{\begin{equation}}
\newcommand{\eeq}{\end{equation}}

\begin{document}
\title[KR Fusion Products]{Fusion products of
  Kirillov-Reshetikhin modules and fermionic multiplicity formulas}

\author{Eddy Ardonne and Rinat Kedem}
\address{EA: Microsoft Station Q, Kavli Institute for Theoretical Physics,
University of California, Santa Barbara, CA 93106 and 
Center for the Physics of Information, California Institute of
Technology, Pasadena, CA 91125, USA. ardonne@kitp.ucsb.edu}
 \address{RK: Department of
     Mathematics, University of Illinois, 1409 W. Green Street, Urbana,
     IL 61801, USA. rinat@uiuc.edu}
\begin{abstract}
  We give a complete description of the graded multiplicity space
  which appears in the Feigin-Loktev fusion product \cite{FL} of
  graded Kirillov-Reshetikhin modules for all simple Lie algebras.
  This construction is used to obtain an upper bound formula for the
  fusion coefficients in these cases.  The formula generalizes the
  case of $\g=A_r$ \cite{AKS2}, where the multiplicities are
  generalized Kostka polynomials \cite{SchWar,KirShi}. In the case of
  other Lie algebras, the formula is the the fermionic side of the
  $X=M$ conjecture \cite{HKOTY}.  In the cases where the
  Kirillov-Reshetikhin conjecture, regarding the decomposition formula
  for tensor products of KR-modules, has been been proven in its
  original, restricted form, our result provides a proof of the
  conjectures of Feigin and Loktev regarding the fusion product
  multiplicites.
\end{abstract}

\maketitle

\section{Introduction}

The Feigin-Loktev fusion product  \cite{FL} is a construction
which defines a $\g$-equivariant grading on the tensor product of
$\g$-modules.
Let $V_1,...,V_N$ be finite-dimensional $\g$-modules as well as
modules over the current algebra $\g[t]:=\g\otimes\C[t]$. Given a set of
pairwise distinct complex numbers $\{\zeta_1,...,\zeta_N\}$, we denote
by $V_i(\zeta_i)$ the $\g[t]$-module localized at $\zeta_i$ (see
Section \ref{localmodule} for definition of a localized module).  

The Feigin-Loktev fusion product of the localized modules
$V_i(\zeta_i)$ is a $\g$-module and a graded $\g[t]$ quotient module.
As a $\g$-module, the fusion product is isomorphic to the usual tensor
product of $\g$-modules, in sufficiently well-behaved cases.  As a
graded $\g[t]$-module, the decomposition of its $n$th graded component into
$\g$-modules is given by the graded multiplicities
$M_{\lambda,\{V_i\}}[n]$:
$$
V_1\ast \cdots \ast V_N (\zeta_1,...,\zeta_N) \simeq \underset{n\geq 0}{\oplus}
\underset{\lambda\in P^+}{\oplus} V_\lambda ^{\oplus M_{\lambda,\{V_i\}}[n]},
$$
where $V_\lambda$ are the irreducible $\g$-modules, and the symbol
$\ast$ denotes the graded (fusion) product rather than the usual tensor
product. The definition of this fusion product depends on the choice of
cyclic vectors of $V_i$ and is given in Section \ref{FLfusiondefn}.

There is no {\it a priori} reason that the graded multiplicities
should be independent of the parameters $\zeta_i$. Nevertheless Feigin and
Loktev conjectured \cite{FL} that they are independent of these
localization parameters. They
suggested that there is a relation between the generating
functions 
\begin{equation}\label{generatingfunction}
M_{\lambda,\{V_i\}}(q):=\sum_{n\geq 0} M_{\lambda,\{V_i\}}[n]\ q^n
\end{equation}
and generalized Kostka polynomials \cite{SchWar,KirShi} (in the
cases where a comparison can be made.)

Several papers contain results which imply special cases of the
Feigin-Loktev conjecture, see e.g.
\cite{FJKLM,Ke,ChLo,AKS2,FKL,FoLit}. For example, in \cite{AKS2}, we
introduced a method to compute the generating function
\eqref{generatingfunction} in the case of $\g=A_r$, for an arbitrary
sequence of irreducible $A_r$-modules with highest weights which are
integer multiples of the fundamental weights. This provided a proof of
the conjectures of Feigin and Loktev by providing explicit,
localization parameter-independent formulas for the graded multiplicities. We
found that they are indeed related to the generalized Kostka
polynomials.

The purpose of this paper is to generalize the results of \cite{AKS2}
to other simple Lie algebras $\g$. We succeed in applying the
techniques of \cite{AKS2} provided that the modules $V_i$ are
$\g[t]$-modules of the Kirillov-Reshetikhin type \cite{Chari}. We find
that the proof of the Feigin-Loktev conjecture in these cases can be
reduced to the proof of the Kirillov-Reshetikhin conjecture. In
general, we obtain an upper bound formula for the fusion multiplicities.

\smallskip

Kirillov-Reshetikhin (KR) modules are the finite-dimensional $\g$- and
$\g[t]$-modules which can be deformed to Yangian modules.  A KR-module
localized at $\zeta\in \C$ is denoted by ${\KR}_{a\omega_i}(\zeta)$,
where $a\in\Z_+$ and $\omega_i$ is one of the fundamental weights of
$\g$.  See Definition \ref{KRmoddefn} for the definition of these
modules in terms of 
generators and relations.
The graded version of this module, $\overline{\KR}_{a\omega_i}$, is
independent of $\zeta$.  It has the same dimension as
$\KR_{a\omega_i}(\zeta)$ and the same $\g$-module structure. The
modules $\overline{\KR}_{a\omega_i}$ are more commonly referred to as a
Kirillov-Reshetikhin modules in the literature \cite{Chari,CM}.  We
refer to this version as a graded KR-module here.

The Kirillov-Reshetikhin conjecture concerns the decomposition of
KR-modules, or the tensor products of several KR-modules, into irreducible
$\g$-modules.  Let us define the graded multiplicities in the fusion
product:
$$
{\KR}_{\mu_1}(\zeta_1)\ast\cdots\ast
{\KR}_{\mu_N}(\zeta_N) \simeq \underset{n\geq
0}{\oplus}\underset{\lambda\in P^+}{\oplus} \
V_{\lambda}^{\oplus M_{\lambda,\{\mu_p\}}[n]} , 
$$
The generating function \eqref{generatingfunction} for the graded fusion
multiplicities is denoted by $M_{\lambda,\{\mu_p\}}(q)$.

The KR-conjecture \cite{KR,HKOTY} gives the total (ungraded)
multiplicities in the form of a restricted sum over binomial
coefficients: 
\begin{equation}\label{krconjecture} 
{\rm dim}\left({\rm Hom}_\g (V_\lambda, \underset{i}{\otimes}
 \KR_{\mu_i}(\zeta_i))\right) =
M_{\lambda,\{\mu_p\}}(1) =
\sum_{\{m_a^{(i)}\in \Z_+ : P_a^{(i)}\geq 0\}}
\prod_{a,i}
{P_a^{i}+m_a^{(i)}\choose m_a^{(i)}} \end{equation}
(see Conjecture
2.3 for the definition of the symbols on the right hand side).

This conjecture implies the completeness of the Bethe ansatz states in
the generalized inhomogeneous XXX spin chain.  The problem of proving it was
addressed in the original work on the subject \cite{KR}, and
subsequently was investigated in various forms by
\cite{Kleber,Chari,HKOTY,Hernandez,Nakajima,FoLit} among others.  To
date, equation \eqref{krconjecture} has only been proven only in
certain special cases (see Section \ref{sec_krcon} for a complete list).

\smallskip 

Our approach in this paper is to describe explicitly the
space dual to the multiplicity space.  It can be expressed as a space
of rational functions with a natural grading.  We then find explicit
formulas for the upper bound of the graded dimension of this space.
The total dimension is, by definition, equal to the multiplicity in the
tensor product decomposition. Since our formulas coincide with the
right-hand side of \eqref{krconjecture} in the limit $q=1$, the proof
of \eqref{krconjecture} would imply that our upper bound is in fact
realized, and is equal the fusion coefficients.

Since the multiplicities depend only on the highest weights
$\mu_p=a_p\omega_{i_p}$ and the Cartan matrix $C$ of $\g$, and are
independent of the localization parameters $\zeta_p$, we thus have a
proof of the Feigin-Loktev conjecture in the cases where
\eqref{krconjecture} has been proven.  That is, the proof of the
Feigin-Loktev conjecture can be reduced to the proof of the
(restricted form of the) Kirillov-Reshetikhin conjecture. We will also
remark below on the way in which a direct proof of the FL-conjecture
can be used to prove the KR-conjecture.

\smallskip

The techniques used in this paper are a generalization of the
techniques introduced in our previous paper \cite{AKS2} for the case
of $\g=A_r$. In that case, the relations for the KR-modules are
precisely the highest-weight conditions for the irreducible
representations of $A_r$ with a highest weight which is a positive
integral multiple of one of the fundamental weights. Therefore, we
were able to compute precisely the fusion Kostka polynomials
\cite{FJKLM}, and show that they are equal to the generalized Kostka
polynomials \cite{KSS}. For other Lie algebras, KR-modules are
precisely those which are amenable to the same type of analysis we
used in \cite{AKS2}. This is due to the form of the relations which
define them, which involve only generators corresponding to simple
roots.

In \cite{FJKLM}, we presented a definition of ``fusion Kostka
polynomials''. The fusion Kostka polynomials are the graded
multiplicities of the fusion of irreducible, finite-dimensional
$\g$-modules. In \cite{FJKLM} we conjectured that these were related
to the generalized Kostka polynomials of \cite{KSS,HKOTY}. This is
true in the case that $\g=A_r$, and was proven in \cite{AKS2}. In
general there is no known fermionic formula for these coefficients,
and we expect that they are different from the fermionic formulas of
\cite{HKOTY} except in the special cases of fusion products of
KR-modules which are also irreducible as $\g$-modules, also known as
miniscule modules.

\smallskip

The paper is organized as follows. The notation and graded KR-modules
localized at $\zeta$ are introduced in Section 2. In Section 3 we
remind the reader of the notion of the graded tensor product of cyclic
$\g[t]$-modules, that is, the fusion product of \cite{FL}.  In Section
4, we describe the dual space to the current algebra
$U(\n_-[t,t^{-1}])$ where $\n_-\subset\g$ is the nilpotent subalgebra.
This is similar to the ideas introduced in \cite{FS} for simple Lie
algebras, and the proofs are a genralization of our constructions in
\cite{AKS2} in the case of $A_r$.  The main difference is the Serre
relation, which governs the structure of singularities in the dual
space. In Section 5, we describe the subspace of matrix elements dual
to the multiplicity space.  We compute a fermionic formula for the associated
graded space, which is an upper bound on these multiplicities.  We
compare the result at $q=1$ to the KR-conjecture, giving an equality
of coefficients in all cases where the equation \eqref{krconjecture}
has been proven.

\vskip.1in
\noindent{\bf Acknowledgements:} We are indebted to David Hernandez,
Peter Littelmann, Masato Okado and 
Anne Schilling for illuminating remarks. RK is funded by NSF grant
DMS-05-00759. EA thanks the KITP for hospitality.

\section{Notation}

\subsection{Finite-dimensional $\g$-modules}
Let $\g$ be a simple Lie algebra of rank $r$, $\Pi =
\{\alpha_1,...,\alpha_r\}$ the set of simple roots and $\Delta$ the
set of all roots. Let $I_r = \{1,...,r\}$ and define the Cartan matrix
$C$ with entries $C_{i,j} = 2(\alpha_i,\alpha_j)/(\alpha_i,\alpha_i)$,
$i,j \in I_r$. Let $\{\omega_1,...,\omega_r\}$ denote the fundamental
weights and $P^+$ the set of dominant integral weights.  

The Cartan decomposition of $\g$ is given by $\g= \n_- \oplus \h
\oplus \n_+$, the generators of $\n_-$ corresponding to simple roots
are denoted by $f_i,\ i\in I_r$, where $f_i\in\g_{-{\alpha_i}}$.
Similarly, the generators of $\n_+$ corresponding to simple roots are
denoted by $e_i$ and the generators of $\h$ by $h_i=[e_i,f_i]$.

Irreducible, finite-dimensional highest weight $\g$-modules are
parameterized by $\lambda\in P^+$ and denoted by $V_\lambda$.  The
module $V_\lambda$ is generated by the action of $U(\n_-)$ on a
highest weight vector $v_\lambda$, which satisfies the properties
\begin{equation}\label{finite}
e_i v_\lambda = 0, \qquad h v_\lambda = \lambda(h) v_\lambda, (h_i\in\h),
\qquad f_i^{l_i+1} v_\lambda = 0 \ ,
\end{equation}
where $\lambda(h_i) = l_i$, and 
the $l_i$ are determined by $\lambda = \sum_{i=1}^r l_i \omega_i$.
We refer to the last condition in \eqref{finite} as the integrability
condition with respect to $\g$.

\subsection{Current algebra modules}\label{localmodule}
Let $t$ be a formal variable, and define the positive current algebra
$\g[t]=g\otimes \C[t]$. We denote its generators by $x[n]:=x\otimes
t^n$, with $x\in\g$ and $n\in\Z_{\geq 0}$. The subalgebra spanned by
$x[0]$ with $x\in \g$ is obviously isomorphic to $\g$, and hence we write
$\g\subset \g[t]$. Therefore, any $\g[t]$-module is also a $\g$-module.

Consider a complex number $\zeta\in \C P$ and the local variable
$t_\zeta = t-\zeta$ ($\zeta\neq \infty$), and $t_\infty = t^{-1}$.
There is a corresponding current algebra $\g[t_\zeta]$, and its
generators are denoted by 
$$x[n]_\zeta := x\otimes (t_\zeta)^n = x\otimes (t-\zeta)^n. $$
Obviously, $x[n]=x[n]_0$. If $\zeta\neq \infty$ then
$\g[t]=\g[t_\zeta]$ as vector spaces.

Let $\zeta\in \C$. Given any $\g[t_\zeta]$-module $V(\zeta)$,
since $\g[t]$ is isomorphic to $\g[t_\zeta]$, there is a natural
action of $\g[t]$ on $V(\zeta)$ which is given by expanding about
$t_\zeta=0$. For any $v \in V(\zeta)$,
\begin{equation}\label{expansion}
x[n]_0\ v = x\otimes t^n\ v = x\otimes(t_\zeta + \zeta)^n\ v
=
\sum_{j\leq n} {n \choose j} \zeta^{j} x[n-j]_\zeta \ v.
\end{equation}
We use the notation $V(\zeta)$ for both the $\g[t_\zeta]$-module and
the $\g[t]$-module with the action described above. It is called the
$\g[t]$-module localized at $\zeta$. Usually, a $\g[t]$-module
localized at $0$ is denoted simply by $V$ instead of $V(0)$.

\subsection{Grading on cyclic $\g[t]$-modules}\label{gradedmodules}
Let $\cA$ be any algebra, then an $\cA$-module is said to by {\em
  cyclic} with cyclic vector $v$ if it is generated by the action of
$U(\cA)$ on $v$.  If $V(\zeta)$ has a cyclic vector $v$ with respect
to $g[t]$ (hence $\g[t_\zeta]$) then we can choose to endow it with a
graded structure with respect to the local variable $t=t_0$, as
follows.

The universal enveloping algebra $U(\g[t])$ is graded by homogeneous
degree in $t$. Although the $\g[t]$-action on $V(\zeta)$ is not graded
if $\zeta\neq 0$, it is filtered. The action of $U(\g[t])$ on
$V(\zeta)$ inherits a filtration from $U(\g[t])$: Let $U^{(\leq n)}$
denote the subspace of $U(\g[t])$ with homogeneous degree in $t$ less
than or equal to $n$.  Define $$\cF^{(n)} = U^{(\leq n)} v.$$
Then
$\cF^{(n)}\subset \cF^{(n+1)}$, $\C v = \cF^{(0)}$ and $V(\zeta) =
\underset{n\geq 0}{\cup} \cF^{(n)}$. Thus, we have defined a filtration of
  $V(\zeta)$ which depends on the choice of cyclic vector.

The associated graded space is
$$
{\rm Gr} V(\zeta) = \underset{n\geq 0}{\oplus} {\rm Gr}[n]
$$
where
$$
{\rm Gr}[n] = \cF^{(n)}/\cF^{(n-1)}.
$$
The structure of the graded space depends on the choice of cyclic
vector $v$. We denote the graded $\g[t]$-module by
$\overline{V}(\zeta)$. It has the structure of a quotient module of
$\g[t]$.

\subsection{Kirillov-Reshetikhin Modules}

The term Kirillov-Reshetikhin (KR) modules originally referred to
certain finite-dimensional modules of the Yangian
\cite{KR}. The Yangian is a deformation of the current algebra
$\g[t]$, where $\g$ is any simple Lie algebra. Thus, one can define
KR-modules for the current algebra as a limiting case of the Yangian
modules.

In \cite{Chari} this idea is used to provide a definition of
KR-modules for $\g[t]$ in terms of current generators and relations.
The result is a finite-dimensional $\g[t]$-module, corresponding to
some highest weight of the form $a\omega_i$ where $\omega_i$ is one of
the fundamental weights of $\g$, and $a\in \Z_+$.  As a $\g$-module, the
KR-module is in fact a highest-weight module with highest-weight
$a\omega_i$.  It is cyclic with respect to the action of $\g[t]$.
However, KR-modules are not necessarily irreducible as $\g$-modules.

We use a definition similar to that found in \cite{Chari}, but we
generalize it to the case where the module is localized at an
arbitrary value of $\zeta\in \C$ to fit the discussion below. (The
definition found in \cite{Chari} corresponds to the special case
$\zeta=1$.)  We also make a distinction between KR-modules and their
associated graded space \cite{CM}, which we call the graded
KR-modules, corresponding to $\zeta=0$ (this is the KR-module defined
in \cite{CM}). As $\g$-modules, these modules
are isomorphic, but they are different as $\g[t]$-modules.

\begin{defn}[Ungraded Kirillov-Reshetikhin modules]\label{KRmoddefn}
Fix $(\zeta,a,i)\in \C^*\times \Z_+\times I_r$. 
Consider the $\g[t]$-module $\KR_{a \omega_i}(\zeta)$
generated by the action of $U(\g[t])$ on the vector $v$ with the properties
\begin{eqnarray}
&& x[n]_\zeta v = 0 \quad {\rm if} \quad n\geq 0,\  x\in \n_+ ;\label{hwcond}\\
&& f_j[n]_\zeta v = 0 \quad {\rm if} \quad n\geq \delta_{i,j}\,
;\label{hwcond2}\\ 
&& f_i[0]_\zeta^{a + 1} v = 0\, ; \\
&& h_j[n]_\zeta v = \delta_{n,0}\delta_{i,j} a v\, . 
\end{eqnarray}
\end{defn}
This implies, for example, that  $f_i[n]_0$
acts on $v$ as
\begin{equation}
f_i[n]_0 v = \zeta^n f_i[0]_\zeta\ v =\zeta^n f_i[0]_0\ v \ .
\end{equation}
Due to the highest weight conditions
\eqref{hwcond}, \eqref{hwcond2} (which imply that $h_j[n]_\zeta v = 0$ if
$n>0$) and the PBW theorem, it is clear that 
 $\KR_{a\omega_i}(\zeta) = U(\n_-[t]) v$.
 
 In the case where $\g=A_r$, $\KR_{a\omega_i}(\zeta) =
 V_{a\omega_i}(\zeta)$, which as an $A_r$-module is the
 irreducible highest-weight module with highest weight $a
 \omega_i$. In general, the $\g[t]$-evaluation module
 $V_{a\omega_i}(\zeta)$, corresponding to the irreducible highest
 weight $\g$-module $V_{a\omega_i}$, is a quotient of the KR-module.
 The weight $a \omega_i$ is the highest $\g$-weight in
 $\KR_{a\omega_i}(\zeta)$.

\begin{defn}[Graded Kirillov-Reshetikhin modules]
  The graded KR-module $\overline{\KR}_{a\omega_i}(\zeta)$
  is the associated graded space
  of the filtered space $\KR_{a\omega_i}(\zeta)$ defined by
  the action of $U(\g[t])$ on the cyclic vector $v$, as in Section
  \ref{gradedmodules}.
\end{defn}

\subsection{The Kirillov-Reshetikhin conjecture}
\label{sec_krcon}

The Kirillov-Reshetikhin conjecture \cite{KR} is an explicit formula
for the decomposition of the tensor product of KR-modules into
irreducible $\g$-modules. (In the literature, there are also other
conjectures which are called the KR-conjecture. We will describe one
of them below.)

Let $\bR$ denote a collection of dominant integral weights of the form
$\{a_p\omega_{i_p}: 1\leq p \leq N\}$, and let $\{\zeta_p:\ 1\leq p\leq N\}$
be a collection of distinct non-zero complex numbers. Define $M_{\lambda, \bR}$
to be the multiplicity of the irreducible $\g$-module $V_\lambda$ in
the tensor product of $N$ KR-modules corresponding to the weights in
$\bR$, localized at $\zeta_p$:
\begin{equation}
M_{\lambda,\bR}=\dim \left({\rm Hom}_\g
\left(\overset{N}
{\underset{p=1}{\otimes}}\KR_{a_p\omega_{i_p}}(\zeta_p),\ V_\lambda\right)  \right)
\end{equation}

Given $\bR$, define $n_a^{(i)}$ to be the number of elements of $\bR$
equal to $a\omega_i$ and define $l^{(i)}$ via $\lambda = \sum_{i=1}^r
l^{(i)} \omega_i$.

\begin{conj}[The Kirillov-Reshetikhin conjecture]
\begin{equation}\label{krconj}
M_{\lambda,\bR} = \sum_{\{m_a^{(i)}\in \Z_+ : P_a^{(i)}\geq 0\}} 
\prod_{a,i}
\binom{P_a^{(i)}+m_a^{(i)}}{m_a^{(i)}}
\end{equation}
where the sum is taken over $m_a^{(i)}$ such that 
$$
\sum_{a} a m^{(i)}_a = \sum_{j,a} C^{-1}_{i,j} a n^{(j)}_a -\sum_j C^{-1}_{i,j} l^{(j)} \ , 
$$
and the ``vacancy numbers'' $P_a^{(i)}$ are
$$
P_a^{(i)} = \sum_{b}\min(a,b)n_b^{(i)} + \sum_{b,j\neq i}\min(|C_{i,j}|
b,|C_{j,i}|a)m_b^{(j)} - 2 \sum_b \min(a,b)  m_b^{(i)},
$$
with $C_{i,j}$ being the Cartan matrix.
\end{conj}

Here, 
\begin{equation}\label{mbinomial}
\binom{n}{m}=\frac{\Gamma(n+1)}{\Gamma(m+1)\Gamma(n-m+1)}.
\end{equation}

This conjecture first appeared in \cite{KR} for the classical Lie
algebras. It has been proven in the following cases: 
\begin{itemize}
\item When $\g=A_r$ \cite{KSS}, and the tensor product of any
  $\KR$-modules. 
\item When $\g=D_r$ and $\bR$ consists of fundamental weights only
\cite{sch}.
\item When $\g$ is any nonexceptional simple Lie algebra 
and $\bR=\{(a_p \omega_1): a_p\in \N\}$ \cite{oss,ss06}.
\end{itemize}
See the review \cite{Schilling} for a full status report.  In each
case, the proof involves a bijection between crystal bases and rigged
configurations \cite{KKR}. We call this version of the
Kirillov-Reshetikhin conjecture KR1.

\begin{remark}
The following is also called the KR-conjecture
\cite{HKOTY,KNT,Nakajima,Hernandez}. It is a modified form of
\eqref{krconj} where the sum is over $P_a^{(i)}\in \Z$, that is, it is
not restricted to be over positive ``vacancy numbers.'' This means
that the sum may contain some negative summands.

We call this second form of the conjecture KR2.  The unrestricted form
of the conjecture can be proven for any Lie algebra $\g$ and any set
of weights \cite{Hernandez}, using an induction system known as the
$Q$-system \cite{Kirillov,Nakajima,KNT}.  On the basis of numerical
evidence, \cite{HKOTY} conjectured that the two forms KR1 and KR2 are
equal in general. There is no direct proof of this fact and it is a
highly nontrivial result.
\end{remark}

\subsection{Affine algebra modules}
The localization procedure of the previous subsection can be extended
to the affine algebra associated with $\g$. For each $\zeta\in\C P$
there is an inclusion
$$
\g[t_\zeta]\subset \g[t_\zeta^{-1},t_\zeta]\subset \g[t^{-1},t]]
$$
where $\C[t^{-1},t]]$ is the space of Laurent series in $t$. 

The algebra $\widetilde{\g}_\zeta=\g[t_\zeta,t_\zeta^{-1}]$ has a
canonical central extension by the central $c$ with the cocycle defined by
\begin{equation}\label{cocycle}
\langle x\otimes f(t), y\otimes g(t)\rangle = \langle x,y\rangle
\oint_{t=\zeta} f'(t) g(t)\ dt, \quad x,y\in\g,\quad f,g\in
\C[t_\zeta,t_\zeta^{-1}].
\end{equation}
We denote the central extention of $\widetilde{\g}_\zeta$ as
$\widehat{\g}_\zeta$. The usual affine algebra is
$\widehat{\g}=\widehat{\g}_0$.
Note that this central extension can be lifted to the algebra $\g[t^{-1},t]]$.

The Cartan decomposition of the affine algebra is $\g = \widehat{\n}_+
\oplus \widehat{\h} \oplus \widehat{\n}_-$ with $\widehat{\n}_{\pm} =
\n_{\pm} \oplus (\g\otimes t^{\pm1}\C[t^{\pm1}])$ and $\widehat{\h} =
\C c \oplus \h$.  

In this paper, we consider only $\widehat{\g}$-modules which are
integrable. Such modules have the property of complete reducibility.
Isomorphism classes of irreducible,
integrable, highest weight modules of $\widehat{\g}$ are indexed by a
positive integer $k$ and certain dominant integral weights $\lambda$ of
$\g$. The integer $k$ is called the {\em level} of the representation:
it is the value of the central element $c$ acting on the irreducible
module. The $\g$-weight $\lambda$ is in the subset $P_k^+\subset P^+$
of weights which have the property $\lambda(\theta) =
\sum_{i=1}^r l_i a_i^\vee \leq k$, where $\theta$ is the highest root
and the $a^\vee_i$ are the co-marks.

Let $V$ be a $\g[t]$-module and let $\widehat{V}$ be a
$\widehat{\g}$-module such that $\widehat{V} = U(\widehat{\n}_-) V$.
(The module $V$ is a {\em top component} of $\widehat{V}$, although
this definition is not unique). In the case of irreducible integrable
modules, the top component is usually taken to be the finite-dimensional
$\g$-module $V_\lambda$, on which $\g[t]$ acts by the evaluation
representation.

A $\widehat{\g}$-module can be localized at $\zeta$ in a similar way
as a $\g[t]$-module, but one must be careful about the appropriate
completion of the algebra. We consider modules of $\widehat{\g}_\zeta$
which are also $\g[t^{-1},t]]$-modules (highest-weight modules have this
property). Currents in $\g[t^{-1},t]]$ act on
$\widehat{\g}_\zeta$-modules by expansion about $t=\zeta$: The formula
\eqref{expansion} can be used even when $n<0$. When acting on any
vector in a highest weight module, the infinite sums are truncated.
The action of the centrally extended current algebra on the module has
a cocycle as in \eqref{cocycle}.  We call the corresponding module the
affine algebra module localized at $\zeta$, denoted by
$\widehat{V}(\zeta)$. If $V(\zeta)$ is a $\g[t]$-module localized at
$\zeta$ then $\widehat{V}(\zeta)$ is the affine algebra module induced
from it by the action of $\g[t^{-1},t]]$.

We refer the reader to the Appendix in \cite{FKLMM} for a further
discussion of tensor (fusion) products of localized affine algebra
modules. In this paper we will use only $\g[t]$-modules
localized at $\zeta$ and affine algebra modules localized at infinity.

We will frequently use generating functions for elements
in $\g[t^{-1},t]]$:
\begin{equation}\label{current}
x(z) = \sum_{n\in\Z} x[n]_0 z^{-n-1}.
\end{equation}
For $\zeta\in \C$ we have
$$
x[n]_\zeta = \sum_j x[n-j]_0 (-\zeta)^j \binom{n}{j} =
\frac{1}{2\pi i} \oint_{z=\zeta} (z-\zeta)^n x(z) dz, \\
$$
where the contour is taken counter-clockwise around the point $\zeta$.
If $\zeta=\infty$, then 
$$
x[n]_\infty = \frac{1}{2\pi i} \oint_{z=\infty} z^{-n} x(z) dz \\
$$
with a clockwise contour.

\section{The fusion product of $\g[t]$-modules and matrix elements}

The idea of the fusion product comes from the fusion product of
integrable $\widehat{\g}$-modules in conformal field theory
\cite{BPZ}. It was reformulated as a graded tensor product of
finite-dimensional $\g[t]$-modules in \cite{FL}. Strictly speaking,
the latter is the restriction of the \cite{BPZ} fusion product to the
top component, or product of primary fields, in the case where the level
is sufficiently large. If the last
condition is not met, we have the level-restricted fusion product,
which will be discussed in a future publication.

\subsection{The Feigin-Loktev graded tensor product}\label{FLfusiondefn}
Let us recall the definition of the \cite{FL}-fusion product,
specialized to the current context. Let $\{V_i(\zeta_i)\}_{i=1}^N$ be
finite-dimensional, cyclic $\g[t]$-modules localized at distinct
points $\zeta_i\in \C$ and let $v_i$ be the cyclic vectors of the
respective modules. Then the tensor product $V_1(\zeta_1)\otimes
\cdots\otimes V_N(\zeta_N)$ is generated by the action of $U(\g[t])$
on the tensor product of cyclic vectors $v_1\otimes \cdots \otimes
v_N$, where the action is by the usual coproduct on the product of
localized modules.

If $v_i$ can be chosen to be highest weight vectors with respect to the
$\g$-action, then the tensor product is generated by the action of the
subalgebra $U(\n_-[t])$ on the tensor product of highest weight
vectors. 

Since the tensor product of modules localized at distinct points is a
cyclic $\g[t]$-module, it has a filtration inherited from $U(\g[t])$,
in the same way as the modules described in Section
\ref{gradedmodules}.  The associated graded space is the graded tensor
product, or the fusion product, of \cite{FL}:
\begin{equation}\label{fusionproduct}
{V}_1 \ast \cdots \ast {V}_N (\zeta_1,...,\zeta_N):=
  {\rm Gr}\left( 
  U(\g[t]) v_1\otimes \cdots 
\otimes v_N\right) .
\end{equation}
Again the filtration, and hence the grading, is $\g$-equivariant, and
therefore the graded components are $\g$-modules. 

Feigin and Loktev conjectured that in certain cases, the graded tensor
product should be independent of the complex numbers $\zeta_i$. Our
purpose in this paper is to prove this conjecture, in the case where
$\g$ is a simple Lie algebra and the modules $V_i$ are KR-modules with
cyclic vectors which are $\g$-highest weight vectors. 

Our method is to compute explicit formulas for the graded fusion
multiplicities. Let
$V[n]$ be the graded component in the graded tensor product:
\begin{equation}
{V}_1 \ast \cdots \ast {V}_N \simeq \underset{n\geq
  0}{\oplus}  V[n],
\end{equation}
where $V[n]$ are finite-dimensional $\g$-modules. Define
\begin{equation}
V[n] \simeq \underset{\lambda\in P^+}{\oplus}
V_\lambda^{\oplus{M_{\lambda,\{V_i\}}[n] }} 
\end{equation}
as $\g$-modules. The generating function for the multiplicities
\begin{equation}
M_{\lambda,\{V_i\}}(q) = \sum_{n\geq 0} M_{\lambda,\{V_i\}}[n] q^n
\end{equation}
is called the graded multiplicity of $V_\lambda$ in the fusion
product.  In \cite{FJKLM}, if $V_i(\zeta_i)$ are taken to be
finite-dimensional irreducible $\g$-modules, the graded multiplicity
was called the fusion Kostka polynomial.

Kirillov-Reshetikhin modules fit the criteria above for the fusion
product: they are finite-dimensional cyclic $\g[t]$-modules. Moreover
their cyclic vector is a highest weight vector with respect to $\g$. Choose
$\bR$ to be a collection of dominant integral highest weights as
before, and choose $\{\zeta_1,...,\zeta_N\}$ to be distinct non-zero
complex numbers.  Define $M_{\lambda,\bR}(q)$ to be the graded
multiplicity of the irreducible $\g$-module $V_\lambda$ in the fusion
product of KR-modules:
$$
M_{\lambda,\bR}(q) = M_{\lambda,\{\KR_{a_p \omega_{i_p}}\}}(q).
$$
Our goal in this paper is to compute the polynomials $M_{\lambda,\bR}(q)$.

\subsection{Fusion products and matrix elements}
It is useful to understand the relation of the graded tensor product
\eqref{fusionproduct} 
of finite-dimensional $\g[t]$-modules to the fusion product of
$\widehat{\g}$-modules: The product \eqref{fusionproduct} is the graded ``top
component'', or the space of conformal blocks,
of the fusion product of $\widehat{\g}$-modules, at sufficiently large
level $k$, together with a grading.  This point of view is taken below
in order to compute the character. See also the Appendix of
\cite{FKLMM} and the introduction in \cite{FJKLM}.

Let $V(\zeta)$ be a graded finite-dimensional cyclic $\g[t]$-module
localized at $\zeta$, and let
$\hV(\zeta)$ be the $\widehat{\g}$-module induced from it at
level $k$. Note
that $V(\zeta)$ need not be an irreducible $\g$-module. Given a collection of
distinct complex numbers $\{\zeta_1,...,\zeta_N\}$ and a collection of
finite-dimensional, cyclic $\g[t]$-modules
$\{V_1(\zeta_1),...,V_N(\zeta_N)\}$, consider the set of induced modules
$\{\hV_1(\zeta_1),...,\hV_N(\zeta_N)\}$. Their fusion product is
denoted by (see Appendix of \cite{FKLMM})
\begin{equation}
\hV_1(\zeta_1)\boxtimes\cdots\boxtimes \hV_N(\zeta_N).
\end{equation}
This is a $\g\otimes \cA_\bzeta$-module, where $\cA_\bzeta$ is the
space of rational functions in $t$ with possible poles at $\zeta_p$
for each $p$. Given a vector $w=w_1\otimes \cdots \otimes w_N$, an
element $x\otimes g(t)$ ($x\in\g, g(t)\in \cA_\bzeta$) acts on the
tensor product by the usual coproduct formula, but the action on the
$p$-th component is given by the expansion of $g(t)$ about $\zeta_p$.

The algebra $\g\otimes\cA_\bzeta\subset \g\otimes\C[t^{-1},t]]$ has a central
extension, with the cocycle given by the sum of the residues at each
point. The level of the fusion product is $k$ if the level of each
module $\hV_p(\zeta_p)$ is $k$. Its decomposition into irreducible
level-$k$-modules of $\widehat{\g}$ is given by the fusion
coefficients or Verlinde numbers. If $k$ is sufficiently large, these
coefficients are equal to the multiplicities of the irreducible
$\g$-modules  in the tensor product of the top components
$V_p(\zeta)$:
\begin{equation}
\label{fusionone}
\mult_{\hV_\lambda(0)}(\hV_1(\zeta_1)\boxtimes\cdots\boxtimes
\hV_N(\zeta_N)) =
\mult_{V_\lambda}(V_1(\zeta_1)\otimes\cdots\otimes V_N(\zeta_N)),\ k \gg 0.
\end{equation}

\begin{remark}
  The only role which the level plays here is the following. It is
  known that any integrable $\widehat{\g}$-module is completely
  reducible, and we use this fact here to argue that the fusion
  product is completely reducible. Hence, it is isomorphic to a direct
  sum of irreducible $\widehat{\g}$-modules localized at $\zeta=0$.
  In fact, we can take $k\to\infty$ in all of the calculations below.
\end{remark}

Thus, the coefficient of the $\g$-module $V_\lambda$ in the decomposition is
the dimension of the space ${\rm Hom}_\g(V_\lambda,
V_1(\zeta_1)\otimes\cdots\otimes V_N(\zeta_N))$. Suppose $V_p$
is generated 
by a highest weight vector $v_p$ for each $p$. Then in terms of matrix
elements, since 
$$V_1(\zeta_1) \otimes \cdots\otimes V_N(\zeta_N) =
U(\n_-[t])v_1\otimes
\cdots \otimes v_N,$$ 
the coefficient of $V_\lambda$ is given by the
inner product with the lowest weight vector $w_{\lambda^*}\in V_\lambda^*$
\begin{equation}\label{matrixelements}
\mult_{V_\lambda}( V_1(\zeta_1)\otimes\cdots\otimes V_N(\zeta_N)) =
\dim \left\{\langle
w_{\lambda^*}, U(\n_-[t]) v_1\otimes \cdots \otimes v_N\rangle\right\}.
\end{equation}

The graded multiplicity of $V_\lambda$ in graded tensor product
is equal to the associated graded space of this space of matrix
elements, where the highest weight vectors are taken to have degree 0
and the grading is inherited from the degree in $t$. 

We may consider $w_{\lambda^*}$ to be a lowest $\g$-weight vector in
the integrable $\widehat{\g}$-module dual to $\widehat{V}_\lambda(0)$,
that is, the module localized at infinity
$\widehat{V}_{\lambda^*}(\infty)$. This module is a graded module with
a grading compatible with the degree in $t$. The graded fusion
coefficients indicate which graded component in a certain
finite-dimensional subspace of 
$\widehat{V}_{\lambda^*}(\infty)$ contain the $\g$-module
$V_\lambda$. One can interpret them as truncated string functions.

\section{The dual space to the algebra}
In order to explain the structure of the space of matrix elements of
the form \eqref{matrixelements}, we will use generating functions for
elements in this space. This is a space of rational functions, with a
certain structure and certain residue properties. In order to explain
this structure we proceed in two steps. In this section, we describe
the dual space to the universal enveloping algebra of
$\n_-[t^{-1},t]$. In the next section we will then restrict this
space to the subspace of matrix elements \eqref{matrixelements}.

\smallskip

Let $\g$ be a simple Lie algebra with Cartan matrix $C$ with entries
given by
$C_{\alpha,\beta}=2(\alpha,\beta)/(\alpha,\alpha)$. The universal
enveloping algebra $U=U(\widetilde{\n}_-)=U(\n_-[t,t^{-1}])$ is
generated by the coefficients of products of generating currents
$f_\alpha(x)$ ($\alpha\in\Pi$). These generating currents satisfy two
types of operator product expansion (OPE) relations. The first comes
from the commutation relation between currents:
\begin{equation}
f_\alpha(x_1) f_\beta(x_2) =
\frac{f_{\alpha+\beta}(x_1)}{(x_1-x_2)}
+\hbox{regular terms}\qquad\hbox{if $\alpha+\beta\in
\Delta$}\label{OPE}
\end{equation}
implying, in particular, that $[f_\alpha(x_1),f_\alpha(x_2)]=0$. The
second relation is a result of the Serre relation for $\g$:
\begin{equation*}
\ad(f_\alpha)^{m_{\alpha,\beta}} f_\beta=0,\quad \alpha,\beta\in \Pi
\end{equation*}
where $m_{\alpha,\beta}=1-C_{\alpha,\beta}$. In terms of currents,
this means that
\begin{equation}
\prod_{i=1}^{m_{\alpha,\beta}}
(x^{(\alpha)}_i-x^{(\beta)})\ f_\alpha(x^{(\alpha)}_1)\cdots
f_\alpha(x^{(\alpha)}_{m_{\alpha,\beta}})
  f_\beta(x^{(\beta)})\Big{|}_{x^{(\alpha)}_1=\cdots=x^{(\alpha)}_{m_{\alpha, \beta}}=x^{(\beta)}}
  =0 \ .
\label{serre}
\end{equation}                                                      

Since $U$ is a graded space with respect to the Cartan subalgebra
$\h\in\g$, define $U[\mathbf m]$ be the subspace consisting of
homogeneous components of weight $-C\mathbf m$ (in the basis of
fundamental weights). Here, $\bm =
(m^{(\alpha_1)},...,m^{(\alpha_r)}).$ That is, $U[\bm]$ 
consists of sums of products of $m^{(\alpha)}$ of each of the currents
$f_\alpha(x^{(\alpha)})$, where $\alpha\in\Pi$.

We want to compute the generating function for the $d$-graded
components of $U[\mathbf m]$, where $d=t \tfrac{d}{dt}$ is the homogeneous
degree in $t$. To do this, we introduce the dual space of the algebra,
and then introduce a filtration on it. This will be the standard
technique for computing characters in everything that follows, and was
described in detail for the case of $\sl_{r+1}$ in \cite{AKS2}. It follows
the general idea described in \cite{FS} for simple $\g$.

The following is a description of the dual space $\mathcal U[\mathbf
  m]$ to $U[\mathbf m]$. It is a space of functions in the set of 
  variables $\{x^{(\alpha)}_i | \alpha\in \Pi, 1\leq i \leq m^{(\alpha)}\}$ with
  the pairing $\mathcal U[\mathbf m]\times U[\mathbf n] \to \C$
  defined as follows. The pairing is $0$ if $\mathbf m \neq \mathbf
  n$, and otherwise it is defined inductively by the relations
\begin{eqnarray}
 \langle 1,1\rangle&=&1\nonumber \\
\langle g, W f_\beta[n]\rangle &=&
\langle\oint_{x^{(\beta)}_{1}=0} dx^{(\beta)}_1
  (x^{(\beta)}_1)^n g(x^{(\beta)}_{1},...), W\rangle,\qquad W\in U[\mathbf
  m-\epsilon_\beta],\nonumber
\end{eqnarray}
where $\epsilon_\beta$ is a standard basis vector. The contour
integral is taken counter-clockwise around the origin in a contour
excluding the singularities at $x^{(\alpha)}_i$ with
$(\alpha,i)\neq(\beta,1)$. Similarly,
\begin{equation*}
\langle g,  f_\beta[n] W\rangle =
\langle\oint_{x^{(\beta)}_1=\infty} dx^{(\beta)}_1
  (x^{(\beta)}_1)^n g(x^{(\beta)}_{1},...), W\rangle,\qquad W\in U[\mathbf
  m-\epsilon_\beta],\nonumber
\end{equation*}
where this time the contour is taken clockwise around the point at
infinity, excluding all the other points $x^{(\alpha)}_i$.

With this pairing, the OPE relation (\ref{OPE}) shows that functions
in the dual space $\mathcal U[\mathbf m]$ may have a simple pole
whenever $x^{(\alpha)}_{i}=x^{(\beta)}_{j}$ if $C_{\alpha,\beta}<0$.
Thus, the dual space is a subspace of 
the space of rational functions of the
form (we impose the usual ordering on the roots)
\begin{equation}
\label{preimage}
g(\bx) =\frac{g_1(\bx)}
{\displaystyle\prod_{\substack{\alpha<\beta\\C_{\alpha,\beta}<0}}\prod_{i,j}
  (x_i^{(\alpha)}-x_j^{(\beta)})}. 
\end{equation}
The function $g_1(\bx)$ is a Laurent polynomial in each of
the variables, symmetric with respect to the exchange of variables of
the same root type $\alpha$, as a result of the fact that generators
of the same root type commute.

Moreover, due to the Serre-type relation (\ref{serre}), the function
$g_1(\bx)$ satisfies the following vanishing condition:
\begin{equation}
\label{serrezeroes}
g_1(\bx)\Big{|}_{x_1^{(\alpha)} = \cdots = x_{m_{\alpha,\beta}}^{(\alpha)}
     = x_1^{(\beta)}} = 0  \  .
\end{equation}

Since these are the only relations among the generating currents
corresponding to simple roots in $\n_-$, this completes the description
of the dual space $\mathcal U[\bm]$.

At this point, we would like to explicitly point out the differences between the
general case $\g$ discussed in this paper, and the case $\sl_{r+1}$ which
we described in \cite{AKS2}. The main difference is, obviously, the
cartan matrix describing the algebra. This difference has important
consequences. On the level of the description it changes the pole structure of
the functions \eqref{preimage}, and the form of the vanishing conditions
\eqref{serrezeroes} for the non simply-laced algebras. More importantly,
it changes the properties of the Kirillov-Reshetikhin modules, which are,
in general, not irreducible, even though, for $\sl_{r+1}$, they are. Nevertheless,
the dual space can be described in one, unified way!

\subsection{Filtration of the dual space}\label{Ufilt}

We introduce a filtration on the space of rational
functions $\mathcal U[\mathbf m]$. The description is very similar to
that explained in
e.g. \cite{AKS2} for the case of $\sl_{r+1}$, except for the vanishing conditions
resulting from the Serre relation. 

Due to the proliferation of indices, we denote the simple roots by their
corresponding index in $I_r$ from this point on.
Let $\boldsymbol \mu = (\mu^{(1)},\cdots,\mu^{(r)})$ be a
multipartition of $\mathbf m=(m^{(1)},\cdots,m^{(r)})$, that is,
$\mu^{(\alpha)}$ is a partition of 
$m^{(\alpha)}$. Define
\begin{equation}
m^{(\alpha)}_a = \#\{\mu_i^{(\alpha)}=a\}
\end{equation}
that is, the number of parts of length $a$ in the partition $\mu^{(\alpha)}$.

Pick a tableau $T$ of shape $\mu^{(\alpha)}$ on the letters
$1,\ldots,m^{(\alpha)}$. Let $a(j)$ be the length of the row in $T$ in
which $j$ appears, and $i(j)-1$ be the number of rows above the row in
which $j$ appears of the same length $a(j)$.
Define the evaluation map
\begin{eqnarray}
\varphi_{\bmu}&:&\ \cU[\bm] \to \cH[\bmu]\\
&& x^{(\alpha)}_j \mapsto y^{(\alpha)}_{a(j),i(j)}\nonumber
\end{eqnarray}
extended by linearity,
where $\cH[\bmu]$ is a space of functions in the variables 
\begin{equation}\label{yvars}
\{ y^{(\alpha)}_{a,i} \ | \alpha\in I_r; \ 1\leq a \leq m^{(\alpha)};\
1\leq i \leq m^{(\alpha)}_a\}.
\end{equation}
There
is a bijective correspondence between the pairs $(a,i)$ for fixed
$\alpha$ and the
parts of the partition $\mu^{(\alpha)}$. 
We order the pairs $(a,i)$ accordingly: 
\begin{equation}\label{pairordering}
(a,i)<(b,j) \quad \hbox{ if $a>b$ or if $a=b$ and $i<j$.}
\end{equation}

The vanishing properties \eqref{serrezeroes} and poles \eqref{preimage}
of functions in $\cU[\bm]$ are
independent of the ordering of the variables $x_{i}^{(\alpha)}$ with
the same index $\alpha$, and therefore the vanishing properties and
poles of the functions in the image in $\cH(\bmu)$ under the
evaluation map are independent of the particular tableau $T$
chosen. Note that the evaluation maps preserve the homogeneous degree.

A lexicographical ordering on multipartitions $\bmu\vdash\bm$ is
defined in the usual way: partitions are ordered lexicographically,
and multipartitions are similarly ordered, $\bnu>\bmu$ if there exists
some $\beta\in I_r$ such that
$\nu^{(\alpha)}=\mu^{(\alpha)}$ for all $\alpha<\beta$ and
$\nu^{(\beta)}>\mu^{(\beta)}$. Let
\begin{equation}
\Gamma_\bmu = \underset{\bnu>\bmu}{\cap} \ker\varphi_{\bnu},\qquad
\Gamma'_\bmu = \underset{\bnu\geq\bmu}{\cap} \ker\varphi_{\bnu}\subset
\Gamma_\bmu. 
\end{equation}
Then $\Gamma_\bmu\subset\Gamma_\bnu$ if $\bmu<\bnu$,
$\Gamma'_{\bmu}=\{0\}$ if $\bmu=
((1^{m^{(1)}}),\cdots,(1^{m^{(r)}}))$, and $\Gamma_\bmu = \cU[\bm]$ if 
$\bmu = ((m^{(1)}),\cdots (m^{(r)}))$. Thus we have a filtration $\cF$
on the space $\cU[\bm]$ parametrized by all multipartitions of $\bm$:
$$
\{0\} \subset \Gamma_{\bmu_1}\subset\cdots
\subset \Gamma_{\bmu_t}= \cU[\bm],
$$
where $\bmu_i< \bmu_{i+1}$ for all $i$.  We can describe the
properties of functions in the image of the associated graded space
$$\Gr\ \cF = \underset{\bmu\vdash \bm}{\oplus}
\Gamma_\bmu/\Gamma_\bmu'
$$
under the induced evaluation maps.

\begin{thm}
\label{thmuim}
The image $\widetilde{\cH}[\bmu]$ of functions in $\Gamma_\bmu/\Gamma'_\bmu$
under the evaluation map $\varphi_\bmu$, is the
space of rational functions of the form
\begin{equation}
\label{umimage}
h(\by)=\frac{\displaystyle
\prod_{\alpha\in I_r}\prod_{(a,i)<(b,j)}
(y_{a,i}^{(\alpha)}-y_{b,j}^{(\alpha)})^{2\min(a,b)} } 
{\displaystyle \prod_{\substack{\alpha<\beta\\
C_{\alpha,\beta}<0}}\prod_{i,j,a,b}
(y_{a,i}^{(\alpha)} - y_{b,j}^{(\beta)})^{\min(|C_{\alpha,\beta}|b, |C_{\beta,\alpha}|a)}}\ h_1(\by),
\end{equation}
where $h_1(\by)$ is an arbitrary Laurent polynomial in the variables
(\ref{yvars}), symmetric with respect to the exchange of variables
$y^{(\alpha)}_{a,i}\leftrightarrow y^{(\alpha)}_{a,j}$. The
induced map
$$
\overline{\varphi}_\bmu: \Gamma_\bmu/\Gamma'_\bmu\to \widetilde{\cH}[\bmu]
$$
is an isomorphism of graded vector spaces.
\end{thm}

We provide a sketch of the proof in the Appendix, since it is very
similar to the proof for the given in \cite{AKS2} in the case of
$\widehat{\sl}_{r+1}$ (in the context of integrable modules).

Thus, we have an isomorphism of graded vector spaces,
$$
\Gr \cF \simeq \underset{\bmu\vdash \bm} {\oplus} \widetilde{\cH}[\bmu].
$$
We will use this fact below to compute the character of a subspace of
$\cU[\bm]$. 

\section{The dual space to the fusion product of KR-Modules}

Our goal is to compute the decomposition of the graded tensor product
of several KR-modules. In order to do this, we introduce the space of
matrix elements, or conformal blocks corresponding to the fusion
product of KR-modules.

\subsection{The space of conformal blocks}
Once more, let $\{\zeta_p
\ | 1\leq p \leq N\}$ be distinct, non-zero complex numbers, and let
$\bR$ be a collection of labels (highest weights) for KR-modules
$$
\bR = \{(a_1 \omega_{\alpha_1}),\cdots, (a_N
\omega_{\alpha_N})\},\quad a_p\in\N,\ \alpha_p\in I_r.
$$
For each $p$, let $\KR_p := \KR_{a_p \omega_{\alpha_p}}(\zeta_p)$
be the KR-module localized at the point $\zeta_p$, with highest
$\g$-weight $a_p \omega_{\alpha_p}$.  We denote by $n_a^{(\alpha)}$
the number of KR-modules with a highest weight equal to
$a\omega_\alpha$. Define $n^{(\alpha)} = \sum_a a n_a^{(\alpha)}$, and
$\bn = (n^{(1)},\cdots n^{(r)})$.

The $q$-multiplicities in the decomposition into irreducible
$\g$-modules of the graded tensor product of $N$ KR-modules
are defined as above:
\begin{equation}\label{qmultiplicities}
{\KR}_1\ast\cdots\ast {\KR}_N \simeq \oplus
M_{\lambda,\bR}(q) V_\lambda, 
\end{equation}
where the notation is understood to mean that a term $q^n V_\lambda$
denotes an irreducible $\g$-module $V_\lambda$ appearing in the graded
component of degree $n$.

To find the multiplicity of $V_\lambda$ in the tensor product, we can
take the inner product with a distinguished vector $w$ in the dual
space $V_\lambda^*$ which we know to have multiplicity 1. For our
purposes here, it is sufficient to take the lowest weight vector
$w_{\lambda^*}\in V_\lambda^*$ (see \cite{AKS2}), which is the lowest
weight vector of the module $V_{-\omega_0(\lambda)}$ where $\omega_0$ is
the longest element in the Weyl group. The vector $w_{\lambda^*}$ has
weight $-\lambda$. 

Let us consider the space of matrix elements of
the form
\begin{equation}\label{matelement}
\langle w_{\lambda^*} f_{\alpha_1}[d_1]\cdots f_{\alpha_r}[d_r]
v_1\otimes \cdots \otimes v_N\rangle.
\end{equation}
where $\alpha_j\in I_r$, and $v_p$ is the highest weight vector
of $\KR_p(\zeta_p)$.

It is convenient to consider the space of generating functions
for such matrix elements:
\begin{equation}
\coinv = \left\{
\langle w_{\lambda^*} |f_{\alpha_1}(x^{(\alpha_1)}_1) \cdots
f_{\alpha_r}(x^{(\alpha_r)}_{i_{m^{(\alpha_r)}}}) v_1\otimes \cdots
\otimes v_N\rangle\right\},\label{cf}
\end{equation}
Here, the variables $x_i^{(\alpha)}$ are formal variables, and the set
includes all possible orderings of the currents corresponding to
simple roots. Coefficients of fixed powers of each of the variables
$x_i^{(\alpha)}$ correspond to matrix elements of the form
\eqref{matelement}. The set of all linearly independent functions in
$\coinv$ corresponds to distinct matrix elements.  We are therefore
interested in the graded dimension of the space of such functions,
which give rise to the polynomials $M_{\lambda,\bR}(q)$, see
Theorem \ref{characters}.

In order to understand such matrix elements,
one should consider ${\KR}_p$ to be  the highest component of a
$\widehat{\g}_\zeta$-module, similar to the ``top component'' of an
irreducible $\widehat{\g}_{\zeta_p}$-module sitting at $\zeta_p$. The induced
$\widehat{\g}_{\zeta_p}$-module is a (graded) direct sum of several
irreducible modules. The graded tensor product of $\KR$-modules is then
the top component of the graded fusion product of several integrable
$\widehat{\g}$-modules in the usual sense of fusion product of
integrable modules%
\footnote{In this paper, we do not consider level restriction, so we always
assume that $k$ is sufficiently large. In practice, this means
$k\geq\sum_\alpha n^{(\alpha)} a_\alpha^\vee$ where $a_\alpha^\vee$
are the comarks.  We will consider the more general, ``level-restricted''
case in a later publication.}.

In this sense, it is clear that one should view $w_{\lambda^*}$ not
just as the lowest weight vector of a $\g$-module $V^*_\lambda$, but
as a vector in the representation $V_{k,\lambda^*}(\infty)$, the
irreducible $\widehat{\g}_\infty$-module, which is annihilated by the
right action of
$\n_-[t]_\infty$.  
This module can also be
considered to be a lowest weight right $\widehat{\g}_0$-module with
lowest weight vector $w_{\lambda^*}$. That means that
\begin{equation}
w_{\lambda^*}f_\alpha[n]_0=w_{\lambda^*}f_\alpha[-n]_\infty = 0\quad
{\rm if}\ n\leq 0.
\end{equation}
Recall that for a module localized at infinity, the local variable is
$t_\infty = t^{-1}$.

The space of matrix elements of the form (\ref{cf}) is a space of
functions, and we analyze the structure of this space below.

\subsection{The dual space of functions}
For convenience, we note the formula for the pairing of a
function $g(\bx)$ in 
the dual space $\cU[\bm]$ with a generator localized at $\zeta\neq 0$:
\begin{equation}
\langle g(\bx), W f_\beta[n]_\zeta\rangle =
\langle\oint_{x^{(\beta)}_1=\zeta} dx^{(\beta)}_1
  (x^{(\beta)}_1-\zeta)^n g(\bx), W\rangle,\qquad W\in U[\bm-\epsilon_\beta]
\end{equation}
and
\begin{equation}
\langle g(\bx),  f_\beta[n]_\infty W\rangle =
\langle\oint_{x^{(\beta)}_1=\infty} dx^{(\beta)}_1
  (x^{(\beta)}_1)^{-n} g(\bx), W\rangle,\qquad W\in U[\bm-\epsilon_\beta].
\end{equation}
The second integral excludes all poles but the one at infinity.
These formulas result from the definitions $x[n]_\zeta = x\otimes
(t-\zeta)^n$ and $x[n]_\infty = x\otimes t^{-n}$.

If $\g(\bx)\in \coinv$ then it should be in the orthogonal
complement to all of the relations in the algebra, the left 
ideal of the algebra
acting on the tensor product to the left and the right ideal of the
algebra acting on the module at infinity to the right. These relations
imply the following.

\begin{enumerate}
\item {\bf Zero weight condition.}
Clearly, the correlation function (\ref{cf}) will vanish unless
  the total weight with respect to $\h\subset\g$ is 0. Denoting $\lambda =
  \sum_\alpha{l^{(\alpha)} \omega_\alpha}$ and $\mathbf l =
  (l^{(1)},\cdots,l^{(r)})$, this means that
  $\bm=(m^{(1)},\cdots,m^{(r)})$ is determined by the 
  equation
\begin{equation}
\label{spin}
\bm = C^{-1} \cdot (\bn-\bl) \ .
\end{equation}
Thus, 
$\coinv\subset \cU[\bm]$ with $\bm$ fixed by (\ref{spin}).
\item {\bf Relations in the algebra.}
Since the space $\coinv$ is a subspace of $\cU[\bm]$,
  any function $g(\bx)\in\coinv$ is of the form
\begin{equation}\label{gpoles}
g(\bx) =
\frac{g_1(\bx)}
{\displaystyle\prod_{\underset{C_{\alpha,\beta}<0}{\alpha<\beta}} 
    \prod_{i,j}(x^{(\alpha)}_i-x^{(\beta)}_j)},
\end{equation}
where $g_1(\bx)$ is a Laurent polynomial which satisfies the
Serre-type relations
\begin{equation}
g_1(\bx)\Big|_{x_1^{(\alpha)}=\cdots = x_{m_{\alpha,\beta}}^{(\alpha)}
  = x_1^{(\beta)}}=0,
\end{equation}
and which is symmetric with respect to exchange of variables of the
same type.
\item {\bf Lowest weight condition.}  The vector $w_{\lambda^*}$ is a
  lowest weight vector with respect to $\g$ of the 
  the module sitting at infinity. The condition at infinity is
  that $w_{\lambda^*} f_\alpha[n]_\infty = 0$ if $n\geq 0$, for any
  $\alpha$. This means that
$$
\oint_{x_1^{(\alpha)}=\infty} dx_1^{(\alpha)} (x_1^{(\alpha)})^{-n}
g(\bx) = 0,\qquad n\geq 0.
$$
Changing variables to $x^{-1}=x_1^{(\alpha)}$, we have
$$
\oint_{x=0} dx x^{-2+n} g(x^{-1},...) = 0, \qquad n\geq 0.
$$
thus, if $\deg_x(g(x,...)) = m$, then $m\leq -2$ if $g$ is the
orthogonal complement of the relations. Therefore we have the degree
restriction for $g(\bx)\in\coinv$:
\begin{equation}
\label{degree}
\deg_{x_i^{\alpha}} g(\bx) \leq -2.
\end{equation}
\item {\bf Highest weight conditions.} 
Each vector $v_p$ is a highest weight vector with weight $\mu_p
  = m\omega_{\alpha_p}$ sitting at $\zeta_p$. Thus, $f_\alpha[n]_{\zeta_p}$
  with $n\geq \delta_{\alpha,\alpha_p}$ generate a left ideal which
  acts trivially on the product $v_1\otimes\cdots\otimes v_N$. This
  means that for functions $g(\bz)\in\coinv$,
\begin{equation}
0=\oint_{x_1^{(\alpha)}=\zeta_p} dx_1^{(\alpha)}
(x_1^{(\alpha)}-\zeta_p)^n g(\bx) \ {\rm if} \ n\geq
\delta_{\alpha,\alpha_p}. 
\end{equation}
Here, the contours are taken so that all other poles are excluded.
This shows that the function $g(\bx)$ may have no poles at
$x_1^{(\alpha)}=\zeta_p$ if $\alpha\neq \alpha_p$. It may have a simple
pole if $\alpha=\alpha_p$. In particular, $g(\bx)$ has no pole at
$x_1^{(\alpha)}=0$ since there is no module localized at 0.
\item {\bf Integrability condition.} 
If $\alpha\neq \alpha_p$ then the function $g(\bx)$ may have a pole
at $x_i^{(\alpha)}=\zeta_p$. However the condition that
\begin{equation}
(f_{\alpha_p}[0]_{\zeta_p})^{a_p+1} v_p = 0
\end{equation}
implies that
\begin{equation}
\oint_{x_1^{(\alpha)}=\zeta_p} dx_1^{(\alpha)}\cdots
\oint_{x_{a+1}^{(\alpha)}=\zeta_p} dx_{a+1}^{(\alpha)}
g(\bx) = 0.
\end{equation}
This shows that if we define $g(\bx)\in \coinv$ by (\ref{gpoles}) with
\begin{equation}
g_1(\bx) = \frac{g_2(\bx)}{\displaystyle \prod_\alpha \prod_p
\prod_{i=1}^{m^{(\alpha)}} (x_i^{(\alpha)} -
\zeta_p)^{\delta_{\alpha_p,\alpha}}}, 
\end{equation}
then the function $g_2(\bx)$ is a polynomial which vanishes on the
following diagonal: 
\begin{equation}
g_2(\bx)\Big|_{x_1^{(\alpha)}=\cdots = x_{a_p+1}^{(\alpha)}=\zeta_p} = 0
\end{equation}
for $p$ such that $\alpha = \alpha_p$.

\end{enumerate}

The space $\coinv$ described above has a filtration by homogeneous
degree in $x_i^{(\alpha)}$, where we take the degree of the factors
$(x_i^{(\alpha)}-\zeta_p)$ in the denominator to be $-1$.
This is equivalent to setting $\zeta_p=0$ for all $p$. 

For any graded space $V$, let ${\rm ch}_q\ V$ denote the generating
function of dimensions of the graded components. That is, $\chq V =
\sum_m q^m \dim V[m] $, where $V[m]$ in the $m$th graded component of
the graded vector space $V$.

The polynomial ${\rm ch}_q\ \coinv$ is related to the polynomial
$M_{\lambda,\bR}(q)$ in equation \eqref{qmultiplicities} as follows.
In the definition of the currents $f_\alpha(x_i^{(\alpha)})$, the
coefficient of the generator $f_\alpha[n]$ is 
$(x_i^{(\alpha)})^{-n-1}$. Therefore,
\begin{equation}
M_{\lambda,\bR} (q) = q^{-|\bm|} \ch_{q^{-1}} \coinv
\end{equation}
Here, $|\bm| = \sum_\alpha m^{(\alpha)}$ is fixed by equation \eqref{spin}.

\subsection{Filtration of the space of matrix elements}

Consider the filtration $\cF$ of the space $\coinv$ defined using the
evaluation mapping, which is inherited from the
filtration of $\cU[\bm]$ introduced in Section \ref{Ufilt}. That is, 
$$\widetilde{\Gamma}_\bmu=\Gamma_\bmu\cap \coinv,\qquad
\widetilde{\Gamma}'_\bmu=\Gamma'_\bmu\cap \coinv.
$$
Functions in the image of the associated graded space of $\cF$
under the evaluation maps $\varphi_\bmu$ can be described explicitly,
in a manner similar to Theorem \ref{thmuim}.
 
\begin{thm}
The image of the space
$\widetilde{\Gamma}_\bmu$ under the evaluation 
map $\varphi_\bmu$ satisfies
$$ 
\varphi_\bmu(\widetilde{\Gamma}_\bmu) \subset
\cH_\cC[\bmu]\subset \cH[\bmu],$$
where $\cH_\cC[\bmu]$ 
is the space of
rational functions $h(\by)$ of the form (\ref{umimage}), where
$h_1(\by)$ is a function of the form
\begin{equation}\label{fusionimage}
h_1(\by)=\prod_p \prod_{a,i}
(y_{a,i}^{(\alpha_p)}-\zeta_p)^{-\min(a_p,a)}\ \bar{h}(\by),
\end{equation}
where $\bar{h}(\by)$ is a polynomial in the variables $\by$,
symmetric under the exchange $y^{(\alpha)}_{a,i} \leftrightarrow
y^{(\alpha)}_{a,j}$, 
such that the total degree of the function $h(\by)$ in the variable
$y_{a,i}^{(\alpha)}$ is less than or equal to $-2a$.
\end{thm}
Again, the proof of this Theorem follows the same arguments as the
proof in the case of $A_r$. In fact, the reason one can describe the
space dual to the fusion product of KR-modules
is that the highest-weight conditions on KR-modules are
identical to the highest-weight conditions on $A_r$-modules with
rectangular highest weights (see equations (5.8) and (5.9) in
\cite{AKS2}). We note that, as in the case of $A_r$, we do not know
how to directly prove the surjectivity of the map
$\varphi_\bmu:\widetilde{\Gamma}_\bmu/\widetilde{\Gamma}_\bmu' \to
\cH_\cC[\bmu]$. Hence, we can only make a statement about the
inclusion of spaces.

\subsection{Character of the fusion product}

The space of functions $\cH_\cC[\bmu]$ is a filtered
space by homogeneous degree in $y_{a,i}^{(\alpha)}$. The associated
graded space is obtained simply by setting $\zeta_p=0$ for all $p$ in
equation \eqref{fusionimage}. That is, it is isomorphic to the space
of functions of the form 
\begin{equation}\label{gradedfunction}
h(\by) = h_0(\by) h_1(\by),
\end{equation} 
with
$$\deg_{y^{(\alpha)}_{a,i}} h(\by)\leq-2a,$$ 
where
\begin{equation}
h_0(\by) = \frac{\displaystyle\prod_\alpha
\prod_{(a,i)<(b,j)}(y_{a,i}^{(\alpha)}- y_{b,j}^{(\alpha)})^{2\min(a,b)}}
{\displaystyle\prod_{p,a,i}(y_{a,i}^{(\alpha_p)})^{\min(a,a_p)}
\prod_{\substack{\alpha<\beta\\C_{\alpha,\beta}<0}}
\prod_{a,i,b,j}(y_{a,i}^{(\alpha)}-y_{b,j}^{(\beta)})^{\min(|C_{\alpha,\beta}|b,|C_{\beta,\alpha}|a)}}.
\end{equation}

Define $$A_{a,b}^{\alpha,\beta}=\delta_{\alpha,\beta}\min(a,b)$$ and
$$B_{a,b}^{\alpha,\beta}=\min(|C_{\alpha,\beta}|b,|C_{\beta,\alpha}|a)
\delta_{\alpha\neq \beta}.$$
We have
$$
\deg_{y^{(\alpha)}_{a,i}}h_0(\by)= -2a-P_a^{(\alpha)}
$$
where
\begin{equation}\label{vacancy}
P_a^{(\alpha)} = [(B-2A)\vec{\mathbf m}+A\vec{\mathbf n}]_a^{(\alpha)}
\end{equation}
where $[\vec{\mathbf m}]_a^{(\alpha)}=m_a^{(\alpha)}$, so that
$[A\vec{\mathbf
m}]_a^{(\alpha)}=\sum_{\beta,b}A_{a,b}^{\alpha,\beta}m_b^{(\beta)}$, etc.. 
Recall that $m_a^{(\alpha)}$ is the number of parts of $\mu^{(\alpha)}$
equal to $a$, whereas $n_a^{(\alpha)}$ is the number of elements in
$\bR$ of the form $(a,\alpha)$. The numbers $P_a^{(\alpha)}$ are
called vacancy numbers in the literature.

{}From this, we see that the function $h(\by)$ in equation
\eqref{gradedfunction} is a polynomial, symmetric in each set of
variables $\{y_{a,i}^{(\alpha)}:\ i=1,...,a\}$, subject to the degree
restriction 
$$
\deg_{y^{(\alpha)}_{a,i}} h(\by)\leq P_a^{(\alpha)}.
$$

Moreover the homogeneous degree of $h_0(\by)$ is
$Q(\vec{\bm},\vec{\bn})-||\vec{\bm}||$, where
$$Q(\vec{\bm},\vec{\bn})=
\vec{\mathbf m}^t A\vec{\mathbf m}-\frac{1}{2} \vec{\mathbf m}^t
B\vec{\mathbf m} - \vec{\mathbf m}^t A\vec{\mathbf n}
$$
and $||\vec{\bm}||=|\bm|=\sum_{a,\alpha}am_a^{(\alpha)}$. Thus, we have
shown that we have the following upper bound on the character of the
fusion product: 
\begin{thm}\label{characters}
\begin{equation}\label{inequality}
M_{\lambda,\bR}(q^{-1})=
q^{|\bm|} \ch_q \coinv \leq 
\sum_{\substack{|\vec{\bm}|=\bm\\ C\bm=(\bn-\bl)}} 
q^{Q(\vec{\bm},\vec{\bn})}
\qbin{\vec{\bm}+\vec{\mathbf P}}{\vec{\bm}}_q,
\end{equation}
where
the $q$-binomial coefficient of vectors is the product
of $q$-binomial coefficients over the entries of the vector,
$$
\qbin{\vec{\bn}}{\vec{\bm}}_q:=\prod_{a,\alpha}\qbin{n_a^{(\alpha)}}
{m_a^{(\alpha)}}_q. 
$$
The notation $|\vec{\bm}|=\bm$ means that $\sum_a a
m_a^{(\alpha)}=m^{(\alpha)}$ for all $\alpha$.
\end{thm}

By an upper bound on a polynomial we mean that each of the
coefficients is bounded from above by the expression on the right. 

Since we have chosen to have no representation localized at 0, the
function $h(\by)$ has no poles at 0. Thus, the $q$-binomial
coefficient which describes the character of the space is the usual
$q$-binomial coefficient,
$$
\qbin{n}{m} = \left\{
\begin{array}{ll}\frac{\displaystyle(q)_n}{\displaystyle (q)_m(q)_{n-m}},&
n\geq m\\ 0 & {\rm otherwise} \end{array}\right.
$$
where $\quad (q)_n =\prod_{i=1}^n(1-q^i)$.

The right hand side of \eqref{inequality} is the $M$-side (that is,
the fermionic sum side) of the $X=M$ conjecture it its original form
\cite{HKOTY}. The fusion coefficients are therefore
bounded by $M(\lambda,W,q^{-1})$ in general, where $W$ is the collection of
KR-modules labelled by $\bR$.

When $q=1$, the right hand side is the fermionic sum of the
Kirillov-Reshetikhin conjecture KR1. The conjecture is that the
multiplicity of the irreducible module $V_\lambda$ in the tensor
product of KR-modules is equal to the right hand side, when
$q=1$. This conjecture has been proven to hold for
several special cases \cite{KSS,sch,oss,ss06}. 
Thus, we have the equality at least in these
cases:
\begin{thm}
The equality holds in equation \eqref{inequality} for the cases of
$\g=A_r$ and $\bR$ consisting of arbitrary KR-modules; for the case of
$\g=D_r$ and $\bR$ consisting of fundamental weights;
and for nonexceptional $\g$ if the weights in $\bR$ are
multiples of the weight $\omega_1$.
\end{thm}

\section{Conclusion}

By studying the space dual to the fusion product of arbitrary
graded Kirillov-Reshetikhin modules of $\g[t]$ for any simple Lie
algebra $\g$, we were able to give an upper bound for the fusion
coefficients.  This upper bound is an equality in the cases where the
Kirillov-Reshetikhin conjecture has been proven in the form KR1.  We
expect that the upper bound provided in this paper is in fact realized
in general. One way to prove this claim is to show that the evaluation
map $\varphi_\bmu$ is surjective. This would also provide an alternative
method of proof of the KR conjecture.

It is clear from the construction we used that the fusion product of
Kirillov-Reshetikhin modules does not depend on the complex numbers
$\zeta_i$ (the `location' of the modules). Thus, in the cases the
Kirillov-Reshetikhin conjecture is proven in the form KR1 (see the
list in section \ref{sec_krcon}), we have proven the conjecture by
Feigin and Loktev, which states that the fusion product of KR-modules
is independent of the $\zeta_i$. 

The results obtained in this paper can be used to provide explicit
character formulas for arbitrary highest weight modules of arbitrary
(non-twisted) affine Lie algebras. In a previous paper \cite{AKS2}, we
showed how this can be done in the case of general $\widehat{\sl}_{r+1}$
modules.  It is possible to generalize the argument of \cite{AKS2} to
arbitrary $\widehat{\g}$.
 We will provide the details in a forthcoming
publication.

In this paper, we did not consider level restriction, that is, we
considered $k\gg 0$. However, the level restricted fusion products of KR-modules
are very interesting. We will study these using methods similar to \cite{AKS2}
in a future publication.

\begin{appendix}
\section{Proof of Theorem \ref{thmuim}}

First, we must prove that the evaluation map ${\varphi}_\bmu:
\Gamma_\bmu\to \widetilde{\cH}[\bmu]$, where the space
$\widetilde{\cH}[\bmu]$ is described in Theorem \ref{thmuim},
is well-defined.

The numerator in equation (\ref{umimage}) describes the vanishing
properties  of
functions in $\varphi_\bmu(\Gamma_\bmu)$.
\begin{lemma}
\label{zerolemma}
Let $g(\bx) \in \Gamma_\bmu$. Then the function $\varphi_\bmu
(g(\bx))$ has a zero of order at least $2 \min(a,b)$ whenever
$y^{(\alpha)}_{a,i}=y^{(\alpha)}_{b,j}$.
\end{lemma} 

\begin{proof}
  The proof of this Lemma is identical to the proof of Lemma (3.7) of
  \cite{AKS2}. Consider the two sets of variables, corresponding to
  two different parts of $\mu^{(\alpha)}$, denoted by $R_a$ and $R_b$,
  of a lengths $a$ and $b$, respectively. Denote the variables
  corresponding to $R_a$ by $\{x_i|\ 1\leq i\leq a\}$ and those
  corresponding to $R_b$ by $\{x'_j|\ 1\leq j\leq b\}$. 
  Without loss of generality, assume $a\geq b$.
  
  Let $\varphi_\bmu = \varphi_2 \circ \varphi_1$, where $\varphi_1$ is
  the evaluation of all the variables except the set corresponding to
  $R_b$. Under $\varphi_1$, the variables $x_i$ are all evaluated at
  $y_a$. The map $\varphi_2$ is the evaluation of the variables $x'_j$
  of $R_b$ at $y_b$.
  
  The image $\bar{g}= \varphi_1 (g(\bx))$ is vanishes whenever
  $y_a = x'_j$ for all $j$, because the evaluation $x'_j=y_a$
  corresponds to an evaluation map $\varphi_{\bmu'}$ for some
  $\bmu'>\bmu$. Functions in $\Gamma_\bmu$ are in the kernel of all
  such evaluation maps, by definition. Thus,
$$
\bar{g} = \prod_{j=1}^{b}({y_a-x'_{j}})\ \bar{g}_1.
$$

Moreover, since $g(\bx)$ is a symmetric under the exchange of
variables corresponding to the same root, it is easy to see that
$\bar{g}_1$ also vanishes whenever $y_a-x'_j$ for each
$j$. It follows that $\varphi_\bmu (g(\bx)) = \varphi_2
(\bar{g})$ has a zero of order at least $2 \min(a,b)$ whenever
$y_a = y_{b}$. The factor $\min(a,b)$ appears because we assumed
that $a\geq b$.
\end{proof}

The denominator in equation (\ref{umimage}) describes the pole
structure of functions in the image of the graded component.
\begin{lemma}
\label{polelemma}
The image under the evaluation map $\varphi_\bmu$ of any function in
$\cU[\bm]$ has a pole of order at most
$\min(|C_{\alpha,\beta}| b, |C_{\beta,\alpha}| a)$ whenever
$y^{(\alpha)}_{a,i} = y^{(\beta)}_{b,j}$
\end{lemma}

\begin{proof}
The order of the pole is limited by the Serre-type relation
(\ref{serrezeroes}). Let $\alpha$ and $\beta$ be two connected roots in
the Dynkin diagram with $|\alpha|\geq |\beta|$. Then $C_{\beta,\alpha}
= -t$ and $C_{\alpha,\beta}=-1$, where $t\in\{1,2,3\}$. The Serre
relation determines the maximal order of the pole of products of the
form 
\begin{equation}\label{fprod}
f_\alpha(z_1) \cdots f_\alpha(z_a) f_\beta(w_1)\cdots f_\beta(w_{b}).
\end{equation}
Note that the order of the pole is independent of the ordering of the
currents in the product.

We assume that $\alpha$ is the longer root, hence the Serre relation
implies that
$$
(z_1-w)(z_2-w) f_\alpha(z_1) f_\alpha(z_2)
f_\beta(w)\Big|_{z_1=z_2=w}=0 
$$
and
$$
(z-w_1)\cdots(z-w_{t+1}) f_\alpha(z)f_\beta(w_1)\cdots
f_\beta(w_{t+1})\Big|_{z=w_1=\cdots=w_{t +1}}=0.
$$

Let $R^{(\alpha)}$ denote a row in $\mu^{(\alpha)}$ of length $a$, and
$R^{(\beta)}$ a row in $\mu^{(\beta)}$ of length $b$. We visualize a pole
as a line connecting box $i$ in $R^{(\alpha)}$ with  box $j$ in $R^{(\beta)}$
(corresponding to the factor $z_i-w_j$ in the denominator). The
allowed configurations of lines connecting boxes are: each box in
$R^{(\beta)}$ has only one line connected to it, and each box in
$R^{(\alpha)}$ may be connected to up to $t$ boxes in $R^{(\beta)}$.

With these allowed configurations, it is easy to see that the maximal
number of lines connecting $R^{(\alpha)}$ with $R^{(\beta)}$ is $\min(b, t
a)$. To see this, first, connect each box in $R^{(\alpha)}$ to one box in
$R^{(\beta)}$. This gives $\min(b,a)$ lines. If $t>1$ and $b>a$ then
connect each box in $R^{(\alpha)}$ to an unconnected box in
$R^{(\beta)}$. This gives $\min(b,2a)$ lines. Finally, if $b>2a$ and
$t=3$, then connect each box in $R^{(\alpha)}$ to an unconnected box in
$R^{(\beta)}$. This gives $\min(b,3a)$ connecting lines.

After evaluating the variables in row $R^{(\alpha)}$ to
$y_{a,i}^{(\alpha)}$ and the variables in $R^{(\beta)}$ to
$y_{b,j}^{(\beta)}$, the maximal pole when these two variables are
equal to each other is of order $\min(|C_{\alpha,\beta}|
b,|C_{\beta,\alpha}| a)$.

\end{proof}

\begin{lemma}
\label{welldef}
The map $\varphi_\bmu: \Gamma_\bmu \rightarrow \widetilde{\cH}[\bmu]$ is well
defined.
\end{lemma}

\begin{proof}
The vanishing properties and pole structure in equation
\eqref{umimage} follow from Lemmas \ref{zerolemma} and
\ref{polelemma}. Furthermore, since functions in $\cU[\bm]$ are
symmetric with respect to exchange of variables with the same root
index $\alpha$, it follows that their evaluation under the map
$\varphi_\bmu$ is symmetric under the exchange of variables
corresponding to the same row of $\mu^{(\alpha)}$.
\end{proof}

This shows that the functions described in Theorem \ref{thmuim} are in
the image  $\varphi_\bmu(\Gamma_\bmu)$. To show that the map $\varphi_\bmu:\
\Gamma_\bmu  \to \widetilde{\cH}[\bmu]$ is surjective, we
give a function in the pre-image of each function in
$\widetilde{\cH}[\bmu]$. The proof of surjectivity is almost identical
to the proof of Lemma 3.16 in \cite{AKS2}, with the appropriate
adjustment made in the pre-image function to account for the different
Serre relations for each $\g$. 

Denote the variables in the preimage of $y_{a,i}^{(\alpha)}$ as follows:
$$
\varphi_\bmu^{-1}(y_{a,i}^{(\alpha)})=
\{x_{a,i}^{(\alpha)}[1],...,x_{a,i}^{(\alpha)}[a]\}\subset
\{x_1^{(\alpha)},\cdots x_{m^{(\alpha)}}^{(\alpha)}\}.
$$
Define
\begin{equation}\label{pre}
g_0(\bx) = \frac{\displaystyle \prod_{\alpha}
\prod_{(a,i)<(a',i')}\prod_{j=1}^{a'}
(x_{a,i}^{(\alpha)}[j]-x_{a',i'}^{(\alpha)}[j])
(x_{a,i}^{(\alpha)}[j+1]-x_{a',i'}^{(\alpha)}[j])}
{\displaystyle \prod_{\substack{\alpha<\beta\\C_{\alpha,\beta}<0}}\prod_{a,i,a',i'}
\prod_{j=1}^{\min(|C_{\alpha,\beta}|b,|C_{\beta,\alpha}|a)}
(x_{a,i}^{(\alpha)}[j]-x_{a',i'}^{(\beta)}[j])}
\end{equation}
where in both the numerator and denominator, we take the argument $j$ to
be a cyclic label: That is, 
$$x_{a,i}^{(\alpha)}[j+a]:=x_{a,i}^{(\alpha)}[j].$$
Recall that according to the ordering \eqref{pairordering},
 if $(a,i)<(a',i')$, then $a\geq a'$.

\begin{thm}
The function $g(\bx)={\rm Sym} (g_0(\bx) g_1(\bx))$, where $g_1(\bx)$ is an
arbitrary Laurent polynomial in the variables $\bx$ and the
symmetrization is over each set of variables with the same label
$\alpha$, is in $\Gamma_\bmu$. Moreover, given a symmetric Laurent
polynomial $h_1(\by)$, there exists a Laurent polynomial $g_1(\bx)$ such
that
\begin{equation}
\varphi_\bmu(g(\bx)) = \xi h(\by)
\end{equation}
where $\xi\in \R$ and $h(\by)$ is the function of equation
 \eqref{umimage}. 
\end{thm}
\begin{proof}
See Lemmas 3.11, 3.13 and 3.15 of \cite{AKS2}. The only difference is
in the vanishing condition due to the Serre relations.
\end{proof}
 
By definition, the induced map
$\overline{\varphi}_\bmu:\Gamma_\bmu/\Gamma'_\bmu$ is injective, since
$\ker(\varphi_\bmu)\cap \Gamma_\bmu = \Gamma'_\bmu$. Thus,
there is an isomorphism:
\begin{lemma}
The induced map
\begin{equation}
\overline{\varphi}_\bmu : \Gr_\bmu \Gamma \rightarrow \widetilde{\cH}[\bmu]
\end{equation}
is an isomorphism of graded vector spaces.
\end{lemma}
This concludes the proof of theorem \ref{thmuim}.

\end{appendix}

\newcommand{\etalchar}[1]{$^{#1}$}


\begin{thebibliography}{HKO{\etalchar{+}}99}

\bibitem[AKS06]{AKS2}
Eddy Ardonne, Rinat Kedem, and Michael Stone.
\newblock Fermionic characters and arbitrary highest-weight integrable affine
  $sl_{r+1}$ modules.
\newblock {\em Comm. Math. Phys.}, 264(2):427--464, 2006.

\bibitem[BPZ84]{BPZ}
A.~A. Belavin, A.~M. Polyakov, and A.~B. Zamolodchikov.
\newblock Infinite conformal symmetry in two-dimensional quantum field theory.
\newblock {\em Nuclear Phys. B}, 241(2):333--380, 1984.

\bibitem[Cha01]{Chari}
Vyjayanthi Chari.
\newblock On the fermionic formula and the {K}irillov-{R}eshetikhin conjecture.
\newblock {\em Internat. Math. Res. Notices}, (12):629--654, 2001.

\bibitem[CL]{ChLo}
Vyjayanthi Chari and Sergei Loktev.
\newblock Weyl, fusion and {D}emazure modules for the current algebra of
  $sl_{r+1}$.
\newblock math.QA/0502165.

\bibitem[CM]{CM}
Vyjayanthi Chari and Adriano Moura.
\newblock The restricted {K}irillov-{R}eshetikhin modules for the current and
  twisted current algebras.
\newblock {\em Comm. Math. Phys.}, 266:431--454, 2006.

\bibitem[FJK{\etalchar{+}}04]{FJKLM}
B.~Feigin, M.~Jimbo, R.~Kedem, S.~Loktev, and T.~Miwa.
\newblock Spaces of coinvariants and fusion product. {I}. {F}rom equivalence
  theorem to {K}ostka polynomials.
\newblock {\em Duke Math. J.}, 125(3):549--588, 2004.

\bibitem[FKL{\etalchar{+}}01]{FKLMM}
B.~Feigin, R.~Kedem, S.~Loktev, T.~Miwa, and E.~Mukhin.
\newblock Combinatorics of the {$\widehat{{sl}}\sb 2$} spaces of coinvariants.
\newblock {\em Transform. Groups}, 6(1):25--52, 2001.

\bibitem[FKL]{FKL}
B. Feigin, A.N. Kirillov and S. Loktev, Combinatorics and geometry of
higher level Weyl modules, math.QA/0503315.

\bibitem[FL99]{FL}
B.~Feigin and S.~Loktev.
\newblock On generalized {K}ostka polynomials and the quantum {V}erlinde rule.
\newblock In {\em Differential topology, infinite-dimensional Lie algebras, and
  applications}, volume 194 of {\em Amer. Math. Soc. Transl. Ser. 2}, pages
  61--79. Amer. Math. Soc., Providence, RI, 1999.

\bibitem[FL05]{FoLit} G. Fourier and P. Littelmann, Weyl modules, affine
  Demazure modules, fusion products and limit constructions, preprint
  math.RT/0509276. 

\bibitem[Her]{Hernandez}
David Hernandez.
\newblock The {K}irillov-{R}eshetikhin conjecture and solutions of $t$-systems.
\newblock math.QA/0501202.

\bibitem[HKO{\etalchar{+}}99]{HKOTY}
G.~Hatayama, A.~Kuniba, M.~Okado, T.~Takagi, and Y.~Yamada.
\newblock Remarks on fermionic formula.
\newblock In {\em Recent developments in quantum affine algebras and related
  topics (Raleigh, NC, 1998)}, volume 248 of {\em Contemp. Math.}, pages
  243--291. Amer. Math. Soc., Providence, RI, 1999.

\bibitem[Ked04]{Ke}
Rinat Kedem.
\newblock Fusion products, cohomology of {${\rm GL}\sb N$} flag manifolds, and
  {K}ostka polynomials.
\newblock {\em Internat. Math. Res. Notices}, (25):1273--1298, 2004.

\bibitem[KKR86]{KKR}
S.~V. Kerov, A.~N. Kirillov, and N.~Yu. Reshetikhin.
\newblock Combinatorics, the {B}ethe ansatz and representations of the
  symmetric group.
\newblock {\em Zap. Nauchn. Sem. Leningrad. Otdel. Mat. Inst. Steklov. (LOMI)},
  155(Differentsialnaya Geometriya, Gruppy Li i Mekh. VIII):50--64, 193, 1986.

\bibitem[Kir83]{Kirillov}
A.~N. Kirillov.
\newblock Combinatorial identities and completeness of states of the
  {H}eisenberg magnet.
\newblock {\em Zap. Nauchn. Sem. Leningrad. Otdel. Mat. Inst. Steklov. (LOMI)},
  131:88--105, 1983.
\newblock Questions in quantum field theory and statistical physics, 4.

\bibitem[KS02]{KirShi} A generalization of the Kostka-Foulkes polynomials.
  J. Algebraic Combin.  15  (2002), 27--69.  

\bibitem[KR90]{KR}
A.N. Kirillov and N.~Yu. Reshetikhin.
\newblock Representations of {Y}angians and multiplicities of the occurrence of
  the irreducible components of the tensor product of representations of simple
  lie algebras.
\newblock {\em J. Soviet Math.}, 52(3):3156--3164, 1990.

\bibitem[KSS02]{KSS}
Anatol~N. Kirillov, Anne Schilling, and Mark Shimozono.
\newblock A bijection between {L}ittlewood-{R}ichardson tableaux and rigged
  configurations.
\newblock {\em Selecta Math. (N.S.)}, 8(1):67--135, 2002.

\bibitem[Kle97]{Kleber}
Michael Kleber.
\newblock Combinatorial structure of finite-dimensional representations of
  {Y}angians: the simply-laced case.
\newblock {\em Internat. Math. Res. Notices}, (4):187--201, 1997.


\bibitem[KNT02]{KNT}
A. Kuniba, T. Nakanishi, and Z. Tsuboi.
\newblock The canonical solutions of the {$Q$}-systems and the
  {K}irillov-{R}eshetikhin conjecture.
\newblock {\em Comm. Math. Phys.}, 227(1):155--190, 2002.


\bibitem[Nak03]{Nakajima}
Hiraku Nakajima.
\newblock {$t$}-analogs of {$q$}-characters of {K}irillov-{R}eshetikhin modules
  of quantum affine algebras.
\newblock {\em Represent. Theory}, 7:259--274 (electronic), 2003.

\bibitem[OSS03]{oss}
Masato Okado, Anne Schilling, and Mark Shimozono.
\newblock A crystal to rigged configuration bijection for nonexceptional affine
  algebras.
\newblock In {\em Algebraic combinatorics and quantum groups}, pages 85--124.
  World Sci. Publishing, River Edge, NJ, 2003.

\bibitem[Sch]{Schilling}
Anne Schilling.
\newblock {$X=M$} theorem: Fermionic formulas and rigged configurations under
  review.
\newblock math.QA/0512161.

\bibitem[Sch05]{sch}
Anne Schilling.
\newblock A bijection between type {$D\sp {(1)}\sb n$} crystals and rigged
  configurations.
\newblock {\em J. Algebra}, 285(1):292--334, 2005.

\bibitem[SF94]{FS}
A.~V. Stoyanovski\u{\i} and B.~L. Fe\u{\i}gin.
\newblock Functional models of the representations of current algebras, and
  semi-infinite {S}chubert cells.
\newblock {\em Funktsional. Anal. i Prilozhen.}, 28(1):68--90, 96, 1994.

\bibitem[SS06]{ss06}
Anne Schilling and Mark Shimozono.
\newblock {$X=M$} for symmetric powers.
\newblock {\em J. Algebra}, 295(2):562--610, 2006.

\bibitem [SW99]{SchWar} Anne Schilling and Ole Warnaar, Inhomogeneous
  lattice paths, generalized Kostka polynomials and A$_{n-1}$
  supernomials, Commun. Math. Phys. 202 (1999) 359-401.
\end{thebibliography}
\end{document}